\documentclass[11pt]{article}
\renewcommand{\baselinestretch}{1.2}
\newcommand{\dated}{\mbox{} \hfill {\small [{\tt \today}]}} \usepackage{amsmath,amssymb,amsfonts,diagrams}
\newenvironment{keywords}{\noindent\small {\it Keywords\/}:}{\vskip 4pt}
\newenvironment{classification}{\noindent\small 2000 {\it Mathematics Subject
Classification\/}:}{\vskip 12pt}

\renewcommand{\iff}{\quad\Longleftrightarrow\quad}
\newcommand{\defiff}{\quad:\Longleftrightarrow\quad}
\newcommand{\comps}{{\mathbb C}}

\newcommand{\ints}{{\mathbb Z}}
\newcommand{\posints}{{\mathbb N}}

\newcommand{\free}{{\mathbb F}}

\newcommand{\tensor}{\otimes}

\newcommand{\Tensor}{\hat{\otimes}}

\newcommand{\cstar}{{C^\ast}}

\newcommand{\id}{{\mathrm{id}}}

\newcommand{\A}{{\mathfrak A}}

\newcommand{\Hilbert}{{\mathfrak H}}

\newcommand{\PM}{\operatorname{PM}}

\newcommand{\supp}{{\operatorname{supp}}}

\newcommand{\cl}[1]{#1^-}

\usepackage{amsthm,enumerate}
\theoremstyle{plain}
\newtheorem{theorem}{Theorem}[section]
\newtheorem{lemma}[theorem]{Lemma}

\newtheorem{proposition}[theorem]{Proposition}
\theoremstyle{definition}
\newtheorem{definition}[theorem]{Definition}
\theoremstyle{remark}
\newtheorem*{remark}{Remark}
\newtheorem*{example}{Example}
\newtheorem*{rems}{Remarks}
\newtheorem*{exs}{Examples}
\newenvironment{remarks}{\begin{rems}\begin{enumerate}}{\end{enumerate}\end{rems}}
\newenvironment{examples}{\begin{exs}\begin{enumerate}}{\end{enumerate}\end{exs}}
\newenvironment{items}{\begin{enumerate}[\rm (i)]}{\end{enumerate}}

\title{Cohen--Host type idempotent theorems \\
for representations on Banach spaces \\
and applications to Fig\`a-Talamanca--Herz algebras}
\author{{\it Volker Runde}\thanks{Research supported by NSERC under grant no.\ 227043-04.}}
\date{}
\begin{document}
\maketitle
\begin{abstract}
Let $G$ be a locally compact group, and let ${\cal R}(G)$ denote the ring of
subsets of $G$ generated by the left cosets of open subsets of $G$. The
Cohen--Host idempotent theorem asserts that a set lies in ${\cal R}(G)$ 
if and only if its indicator function is a coefficient function of a unitary
representation of $G$ on some Hilbert space. We prove related results for 
representations of $G$ on certain Banach spaces. We apply our Cohen--Host
type theorems to the study of the Fig\`a-Talamanca--Herz algebras $A_p(G)$ 
with $p \in (1,\infty)$. For arbitrary $G$, we characterize those closed 
ideals of $A_p(G)$ that have an approximate identity bounded by $1$ in terms of their 
hulls. Furthermore, we characterize those $G$ such that $A_p(G)$ is 
$1$-amenable for some --- and, equivalently, for all --- $p \in
(1,\infty)$: these are precisely the abelian groups.
\end{abstract}
\section*{Introduction}
Let $G$ be a locally compact abelian group with dual group $\hat{G}$. In
\cite{Coh}, P.\ J.\ Cohen characterized the idempotent elements of the
measure algebra $M(G)$ in terms of their Fourier--Stieltjes transforms:
$\mu \in M(G)$ is idempotent if and only if $\hat{\mu}$ is the indicator
function of a set in the coset ring of $\hat{G}$.
\par
In \cite{Eym}, P.\ Eymard introduced, for a general locally compact group
$G$, the Fourier algebra $A(G)$ and the Fourier--Stieltjes algebra $B(G)$.
If $G$ is abelian, the Fourier and Fourier--Stieltjes transform, respectively,
yield isometric Banach algebra isomorphisms $A(G) \cong L^1(\hat{G})$ 
and $B(G) \cong M(\hat{G})$. (In the framework of Kac algebras,
this extends to a duality between $L^1(G)$ and $A(G)$ for arbitrary $G$;
see \cite{ES}). In \cite{Hos}, B.\ Host extended Cohen's idempotent theorem to 
Fourier--Stieltjes algebras of arbitrary locally compact groups. Besides
being more general than Cohen's theorem, Host's result also has a much simpler
proof that only requires elementary operator theory on Hilbert spaces.
\par
Host's result --- to which we shall refer as to the \emph{Cohen--Host 
idempotent theorem} or simply the \emph{Cohen--Host theorem} --- 
has turned out to be a tool of crucial importance
in the investigation of $A(G)$ and $B(G)$. We mention only three recent
applications:
\begin{itemize}
\item \textit{Homomorphisms between Fourier algebras.} Already Cohen used
his theorem to describe the algebra homomorphism from $A(G)$ to $B(H)$,
where $G$ and $H$ are locally compact abelian groups (\cite{Coh2}). 
Host extended Cohen's result to a setting where only $G$ had to be abelian 
(\cite{Hos}). This line of research culminated only recently with a complete
description of the completely bounded algebra homomorphism from $A(G)$ to
$B(H)$ with $G$ amenable (\cite{Ili} and \cite{IS}).
\item \textit{Ideals of $A(G)$ with a bounded approximate identity.} In
\cite{FKLS}, the closed ideals of $A(G)$ (for amenable $G$) that have a 
bounded approximate identity are completely characterized in terms of their
hulls. One direction of this result requires operator space methods
(see \cite{ER}), but the converse implication relies mainly on the
Cohen--Host result.
\item \textit{Amenability of $A(G)$.} B.\ E.\ Forrest and the author,
in \cite{FR}, characterized those locally compact groups $G$ for which
$A(G)$ is amenable in the sense of \cite{Joh1}: they are precisely those
with an abelian subgroup of finite index (\cite[Theorem 2.3]{FR}). The
proof in \cite{FR} relies on the Cohen--Host idempotent theorem only
indirectly --- through \cite{FKLS} ---, but recently, the author gave an
alternative proof that invokes the idempotent theorem directly
(\cite{RunPP}).
\end{itemize}
\par
The Fig\`a-Talamanca--Herz algebra $A_p(G)$ for $p \in (1,\infty)$ were
introduced and first studied by C.\ Herz (\cite{Her1} and \cite{Her2});
more recent papers investigating those algebras are, 
for example, \cite{For2}, \cite{Gra}, \cite{LNR}, and \cite{Mia}.
They are natural generalizations of $A(G)$ in the sense that 
$A_2(G) = A(G)$. The algebras $A_p(G)$ share many properties of $A(G)$. For 
instance, Leptin's theorem
(\cite{Lep}) extends easily to $A_p(G)$ (\cite{Her2}): $G$ is amenable if
and only if $A_p(G)$ has an approximate identity for some --- and, 
equivalently, for all --- $p \in (1,\infty)$. Nevertheless, since $A_p(G)$ has 
no obvious connection with Hilbert space for $p \neq 2$, the powerful methods 
of operator algebras are not available anymore --- or are at least 
not as easily applicable --- for the study of general Fig\`a-Talamanca--Herz 
algebras. As a consequence, many questions to which the 
answers are easy --- or have 
at least long been known --- for $A(G)$ are still open for general $A_p(G)$. 
\par
Since the Cohen--Host theorem is about $B(G)$ its use for the investigation
of Fig\`a-Talamanca--Herz algebras is usually limited to the case where 
$p = 2$. In this paper, we therefore strive for
extensions of this result that are applicable to the study of $A_p(G)$ for
general $p \in (1,\infty)$. The elements of $B(G)$ can be interpreted
as the coefficient functions of the unitary representations of $G$ on 
Hilbert spaces, so that the Cohen--Host theorem (or rather its difficult
direction) can be formulated as follows: if the indicator function of a
subset of $G$ is a coefficient function of a unitary representation of 
$G$ on a Hilbert space, then the set lies in the coset ring of $G$. We shall 
prove two Cohen--Host type theorems for representations on Banach spaces.
In particular, we shall extend \cite[Theorem 2.1]{IS} to isometric 
representations on Banach spaces which are smooth or have a smooth dual. 
\par
We apply this Cohen--Host type theorem to the study of 
general Fig\`a-Talamanca--Herz algebras.
\par
First, we characterize
those closed ideals $I$ of $A_p(G)$ that have an approximate identity bounded
by $1$: we shall see that $I$ has such an approximate identity if and only if
$I$ consists precisely of those functions in $A_p(G)$ that vanish outside
some left coset of an open, amenable subgroup of $G$. This result is related to
\cite[Theorems 2.3 and 4.3]{FKLS}, and extends \cite[Proposition 3.12]{For1}
from $p=2$ to arbitrary $p \in (1,\infty)$.
\par
Secondly, we extend \cite[Theorem 3.5]{RunPP} and show that $A_p(G)$ is 
$1$-amenable for some --- and, equivalently, for all --- $p \in (1,\infty)$ 
if and only if $G$ is abelian.
\subsubsection*{Acknowledgments}
I would like to thank Nico Spronk for suggesting Theorem \ref{monico} to me ---
as well as for detecting a serious gap in an earlier version of this paper --- and Matthew Daws for bringing \cite{BFGM} to my attention. 
\section{Cohen--Host type idempotent theorems for representations on Banach spaces}
Our notion of a representation of a locally compact group on a Banach space
is the usual one:
\begin{definition}
Let $G$ be a locally compact group. Then $(\pi,E)$ is said to be a
\emph{representation} of $G$ on $E$ if $E$ is a Banach space and $\pi \!:
G \to {\cal B}(E)$ is a group homomorphism into the invertible operators
on $E$ which is continuous with respect to the given topology on $G$ and the
strong operator topology on ${\cal B}(E)$. We call $(\pi,E)$ 
\emph{uniformly bounded} if $\sup_{x \in G} \| \pi(x) \| < \infty$ and
\emph{isometric} if $\pi(G)$ only consists of isometries.
\end{definition}
\begin{remarks}
\item Suppose that $(\pi,E)$ is uniformly bounded. Then
\[
  ||| \xi ||| := \sup_{x \in G} \| \pi(x) \xi \| \qquad (\xi \in E) 
\]
defines an equivalent norm on $E$ such that $(\pi,(E,||| \cdot |||))$ is
isometric. This, however, may obscure particular geometric features of the
original norm.
\item Since invertible isometries on a Hilbert space are just the unitary 
operators, the isometric representations of $G$ on a Hilbert space, are just 
the usual unitary representations.
\item Every representation $(\pi,E)$ of $G$ induces a representation of
the group algebra $L^1(G)$ on $E$ through integration, which we denote likewise
by $(\pi,E)$.
\end{remarks}
\par
We are interested in certain functions associated with representations:
\begin{definition}
Let $G$ be a locally compact group, and let $(\pi,E)$ a representation of
$G$. A function $f \!: G \to \comps$ is called a \emph{coefficient function}
of $(\pi,E)$ if there are $\xi \in E$ and $\phi \in E^\ast$ such that
\begin{equation} \label{coeff}
  f(x) = \langle \pi(x) \xi, \phi \rangle \qquad (x \in G).
\end{equation}
If $\| \xi \| = \| \phi \| = 1$, we call $f$ \emph{normalized}.
\end{definition}
\par
The coefficient functions of the unitary representations of a locally compact
group $G$ form an algebra (under the pointwise operations), the 
\emph{Fourier--Stieltjes algebra} $B(G)$ of $G$ (see \cite{Eym}). Moreover, $B(G)$ can be
identified with the dual space of the full group $\cstar$-algebra of $G$, 
which turns it into a commutative Banach algebra.
\par
Extending earlier work by Cohen in the abelian case (\cite{Coh}), Host identified the idempotents of $B(G)$ (\cite{Hos}). Since $B(G)$ consists of continuous
functions, it is clear that an idempotent of $B(G)$ has to be the indicator
function $\chi_C$ of some clopen subset $C$ of $G$. Let ${\cal R}(G)$ denote
the \emph{coset ring} of $G$, i.e.\ the ring of sets generated by all left
cosets of open subgroups of $G$. In \cite{Hos}, Host showed that the 
idempotents of $B(G)$ are precisely of the form $\chi_C$ with 
$C \in {\cal R}(G)$.
\par
Given a representation $(\pi,E)$, where $E$ is not necessarily a Hilbert space,
the set of coefficient functions of $(\pi,E)$ need not be a linear space 
anymore, let alone an algebra. Nevertheless, it makes sense to attempt to
characterize those subsets $C$ of $G$ for which $\chi_C$ is a coefficient
function of $(\pi,E)$.
\par
Without any additional hypotheses, we cannot hope to
extend the Cohen--Host theorem:
\begin{example}
Let $G$ be any locally compact group, and let ${\cal C}_\mathrm{b}(G)$ denote
the bounded, continuous function on $G$. For any function $f \!: G \to \comps$
and $x \in G$, define 
\[
  r_x f \!: G \to \comps, \quad y \mapsto f(yx),
\]
and call $f \in {\cal C}_\mathrm{b}(G)$ \emph{right uniformly continuous} if 
the map
\[
  G \to {\cal C}_\mathrm{b}(G), \quad x \mapsto r_x f
\]
is continuous with respect to the given topology on $G$ and the norm topology
on ${\cal C}_\mathrm{b}(G)$. The set of all right uniformly continuous function
on $G$ is a $\cstar$-subalgebra of ${\cal C}_\mathrm{b}(G)$, which we denote
by $\mathit{RUC}(G)$. Define an isometric representation 
$(\rho,\mathit{RUC}(G))$ by letting $\rho(x)f := r_x f$ for $x \in G$ and
$f \in \mathit{RUC}(G)$. It is then immediate that
\[
  f(x) = \langle \rho(x)f, \delta_e \rangle \qquad (f \in \mathit{RUC}(G), \,
                                                    x \in G),
\]
where $\delta_e$ is the point mass at the identity of $G$, so that every
element of $\mathit{RUC}(G)$ is a coefficient function of 
$(\rho,\mathit{RUC}(G))$. For discrete $G$, it is clear that 
$\mathit{RUC}(G) = \ell^\infty(G)$, so that $\chi_C$ is a coefficient function
of $(\rho,\mathit{RUC}(G))$ for every $C \subset G$.
\end{example}
\par
If we impose restrictions on both the group and the Banach space on which it is
represented, an extension of the Cohen--Host theorem is surprisingly
easy to obtain. 
\par
For any locally compact group $G$, denote the component of the identity
by $G_e$; it is a closed, normal subgroup of $G$. 
Recall that $G$ is said to be \textit{almost connected} if $G/G_e$ is compact.
\begin{theorem} \label{cohost1}
Let $G$ be an almost connected locally compact group.
Then the following are equivalent for $C \subset G$:
\begin{items}
\item $C \in {\cal R}(G)$;
\item $\chi_C \in B(G)$;
\item $\chi_C$ is a coefficient function of a uniformly bounded representation
$(\pi,E)$ of $G$, where $E$ is reflexive.
\end{items} 
\end{theorem}
\begin{proof}
(i) $\Longleftrightarrow$ (ii) is the Cohen--Host
theorem, and (ii) $\Longrightarrow$ (iii) is straightforward.
\par
(iii) $\Longrightarrow$ (i): Let $\xi \in E$ and $\phi \in E^\ast$ such that 
$\chi_C$ is of the form (\ref{coeff}). We can suppose without loss of 
generality that $\{ \pi^\ast(f)\phi : f \in L^1(G) \}$ is dense in $E^\ast$:
otherwise, replace $E^\ast$ by $\cl{\{ \pi^\ast(f)\phi : f \in L^1(G) \}}$ and
$E$ by its quotient modulo $\{ \pi^\ast(f)\phi : f \in L^1(G) \}^\perp$.
\par
We claim that $\mathbb{I} := \{ \pi(x) \xi : x \in G \}$ is uniformly 
discrete in the 
norm topology. To see this, let $x_1,x_2 \in G$ be such that 
$\| \pi(x_1) \xi - \pi(x_2) \xi \| < \frac{1}{C\| \phi \| + 1}$, where
$C := \sup_{x \in G} \| \pi(x) \|$. We thus have
\begin{equation} \label{close}
  | \langle \pi(y)\pi(x_1)\xi, \phi \rangle - 
  \langle \pi(y)\pi(x_2)\xi, \phi \rangle | 
  < 1
  \qquad (y \in G).
\end{equation}
On the other hand, since $\langle \pi(x)\xi,\phi \rangle \in \{ 0,1 \}$ for 
$x \in G$, it is clear that $| \langle \pi(y)\pi(x_1)\xi, \phi \rangle - 
\langle \pi(y)\pi(x_2)\xi, \phi \rangle | \geq 1$ whenever $y \in G$ is
such that $\langle \pi(y)\pi(x_1) \xi, \phi \rangle \neq
\langle \pi(y) \pi(x_2), \phi \rangle$. Combining this with (\ref{close})
yields
\begin{equation} \label{disc1}
  \langle \pi(y)\pi(x_1) \xi, \phi \rangle = 
  \langle \pi(y) \pi(x_2), \phi \rangle
  \qquad (y \in G).
\end{equation}
Integrating (\ref{disc1}) with respect to $y$, we obtain
\begin{multline*}
  \langle \pi(x_1)\xi, \pi(f)^\ast \phi \rangle = 
  \langle \pi(f) \pi(x_1)\xi, \phi \rangle = 
  \langle \pi(f) \pi(x_2)\xi, \phi \rangle =
  \langle \pi(x_2)\xi , \pi(f)^\ast \phi \rangle
  \\ (f \in L^1(G)).
\end{multline*}
Since $\cl{\{ \pi^\ast(f)\phi : f \in L^1(G) \}} = E^\ast$, the Hahn--Banach
theorem yields that $\pi(x_1) \xi = \pi(x_2) \xi$. 
\par
Since $\{ \pi(x) \xi : x \in G_e \}$ is connected in the norm topology of 
$E$, we conclude that $\pi(x)\xi = \xi$ for all $x \in G_e$. 
As a consequence, $\pi(x)\xi$ with $x \in G$ only
depends on the coset of $x$ in $G / G_e$. Hence, the map
\[
  G/G_e \to E, \quad xG_e \mapsto \pi(x)\xi
\]
is well defined, is continuous with respect to the norm topology on $E$,
and clearly has $\mathbb I$ as its range. Since $G/G_e$ is compact, it follows
that $\mathbb I$ is compact and thus finite.
\par
Let $G_d$ denote the group $G$ equipped with the discrete topology.
Define a unitary representation $\tilde{\pi}$ of $G_d$ on $\ell^2(\mathbb{I})$ 
by letting
\[
  \tilde{\pi}(x) \delta_\eta := \delta_{\pi(x) \eta} 
  \qquad (x \in G, \, \eta \in \mathbb{I})
\]
Since $\mathbb I$ is finite, the restriction of $\phi$ to $\mathbb I$
can be identified with an element of $\ell^2(\mathbb{I})^\ast$, which we denote
by $\tilde{\phi}$. By construction, we have
\[
  \left\langle \tilde{\pi}(x) \delta_{\xi}, \tilde{\phi} \right\rangle =
  \left\langle \pi(x) \xi, \phi \right\rangle = \chi_C(x) \qquad (x \in G),
\]
so that $\chi_C \in B(G_d)$. Since $C$ is clopen, $\chi_C$ is continuous, so
that actually $\chi_C \in B(G)$ by \cite[(2.24) Corollaire 1]{Eym}. From
\cite{Hos}, we conclude that $C \in {\cal R}(G)$. 
\end{proof}
\par
In \cite{IS}, M.\ Ilie and N.\ Spronk proved a variant of the Cohen--Host
theorem for normalized coefficient functions in the Fourier--Stieltjes
algebra: they showed that these are precisely the indicator functions of
left cosets of open subgroups (\cite[Theorem 2.1]{IS}). 
As Spronk pointed out to the author, the argument used in \cite{IS} can
be adapted to certain Banach spaces.
\par
The following definition is crucial (see \cite{JL}, for instance):
\begin{definition}
A Banach space $E$ is said to be \emph{smooth} if, for each $\xi \in E
\setminus \{ 0 \}$, there is a unique $\phi \in E^\ast$ such that
$\| \phi \| = 1$ and $\langle \xi, \phi \rangle = \| \xi \|$.
\end{definition}
\par
Extending \cite[Theorem 2.1]{IS}, we obtain:
\begin{theorem} \label{monico}
Let $G$ be a locally compact group. Then the following are equivalent 
for $C \subset G$:
\begin{items}
\item $C$ is a left coset of an open subgroup of $G$;
\item $\chi_C \in B(G)$ with $\| \chi_C \| = 1$;
\item $\chi_C \not\equiv 0$ is a normalized coefficient function of an isometric 
representation $(\pi,E)$ of $G$, where $E$ or $E^\ast$ is smooth.
\end{items}
\end{theorem}
\begin{proof}
(i) $\Longleftrightarrow$ (ii) is \cite[Theorem 2.1(i)]{IS}, 
and (ii) $\Longrightarrow$ (iii) is obvious.
\par
(iii) $\Longrightarrow$ (i): Suppose that $\chi_C \not\equiv 0$ is of the form 
(\ref{coeff}) with $\xi \in E$ and $\phi \in E^\ast$ such that 
$\| x \| = \| \phi \| = 1$. 
\par
We first treat the case where $E^\ast$ is smooth. Fix $x \in C$, and set
\[
  H := \{ y \in G : xy \in C \}.
\]
By definition, we have for $y \in G$ that
\[
  y \in H  \iff \langle \pi(xy)\xi,\phi \rangle =
                \langle \pi(y)\xi, \pi(x)^\ast \phi \rangle =
                \langle \xi, \pi(x)^\ast \phi \rangle = 1. 
\]
Since $E^\ast$ is smooth, there is a \emph{unique} 
$\Psi \in E^{\ast\ast}$ such that $\langle \pi(x)^\ast \phi, \Psi \rangle = 1$.
From this uniqueness assertion, it follows that $\Psi = \xi = \pi(y)\xi$ for
all $y \in H$ and that 
\[
  H = \{ y \in G : \pi(y) \xi = \xi \}.
\]
Consequently, $H$ is a subgroup of $G$, and it is immediate that $C = xH$.
Since $\chi_C$ is continuous, $C$ --- and thus $H$ --- is clopen. This proves
(i).
\par
If $E$ is smooth, an analogous argument yields that there are $x \in G$ and
an open subgroup $H$ of $G$ such that $C = Hx$. Since $H x = x( x^{-1} H x)$,
this also proves (i).
\end{proof}
\par
At this point, we take a look at those spaces to which we shall
apply Theorem \ref{monico} in the next section:
\begin{example}
The \emph{modulus of convexity} of a Banach space $E$ is defined, for
$\epsilon \in (0,2]$ as
\[
  \delta_E(\epsilon) := \inf \left\{ 1 - \frac{\| \xi + \eta \|}{2} : 
  \xi, \eta \in E, \,
  \| \xi \| \leq 1, \, \| \eta \| \leq 1, \| \xi - \eta \| 
  \geq \epsilon \right\} > 0;
\]
if $\delta_E(\epsilon) > 0$ for each $\epsilon \in (0,2]$, then $E$ is called
\emph{uniformly convex} (\cite[Definition 9.1]{Fab}). All uniformly
convex Banach spaces are reflexive (\cite[Theorem 9.12]{Fab}).
If $X$ is any measure space and $p \in (1,\infty)$, then $L^p(X)$
is uniformly convex (\cite[Theorem 9.3]{Fab}). More generally, whenever
$E$ is a uniformly convex Banach space, $X$ is any measure space, and $p
\in (1,\infty)$, the vector valued $L^p$-space $L^p(X,E)$ is again
uniformly convex (\cite{Day}); in particular, for any two measure spaces
$X$ and $Y$ and $p,q \in (1,\infty)$, the Banach space $L^p(X,L^q(Y))$ is
uniformly convex. If $E$ is uniformly convex, then $E^\ast$ is smooth
by (\cite[Lemma 8.4(i) and Theorem 9.10]{Fab}). Hence,
if $G$ is a locally compact group and $C \subset G$ is such that
$\chi_C$ is a normalized coefficient function of an isometric representation
on a Banach space, which is uniformly convex or has a uniformly convex dual, 
then $C$ is a left coset of an open subgroup of $G$ by Theorem \ref{monico}. 
\end{example}
\par
The proof of the general Cohen--Host theorem from \cite{Hos} relies heavily on some (elementary) facts on Hilbert space operators, for which
there is no analog in a more general Banach space setting. Concluding this section, we shall see that there \emph{is} no general Cohen--Host theorem for isometric representations on uniformly convex Banach spaces, even if we demand that the dual spaces be uniformly convex, too:
\begin{example}
A subset $L$ of a group $G$ is called a \emph{Leinert set} (see
\cite{Lei1} and \cite{Lei2}) if, for any $x_1, \ldots, x_{2n} \in L$
with $x_j \neq x_{j+1}$ for $j=1, \ldots, 2n-1$, we have
$x_1^{-1} x_2 x_3^{-1} \cdots x_{2n-1}^{-1} x_{2n} \neq e$. For instance,
the subset $\{ a^n b^n : n \in \ints \}$ of the free group $\free_2$ generated 
by $a$ and $b$ is a Leinert set (\cite[(1.10)]{Lei1}). By the proof of
\cite[(12) Korollar]{Lei2}, the indicator function of an infinite Leinert 
subset of $\free_2$ does not lie in $B(\free_2)$, so that the set does not
belong to ${\cal R}(\free_2)$. On the other hand, the indicator function of
every Leinert subset of a group $G$ is a coefficient function of a uniformly
bounded representation $(\pi,\Hilbert)$ of $G$, where $\Hilbert$ is some Hilbert space (\cite[1.1 Theorem]{Fen}).
By \cite[Proposition 2.3]{BFGM}, there is an equivalent norm $\| \cdot \|$ on $\Hilbert$ such that both $E := (\Hilbert, \| \cdot \|)$ and $E^\ast$ are uniformly convex and such that $(\pi,E)$ is isometric. Consequently, if $L \subset \free_2$ is an infinite Leinert set, then $L \notin {\cal R}(\free_2)$, but $\chi_L$ is a coefficient function of some isometric representation $(\pi,E)$, where both $E$ and $E^\ast$ are uniformly convex Banach spaces.
\end{example}
\section{Applications to $A_p(G)$}
We shall now turn to applications of Theorem \ref{monico} 
to Fig\`a-Talamanca--Herz algebras on locally compact groups.
\par
Let $G$ be a locally compact group. For any function $f \!: G \to \comps$, we 
define $\check{f} \!: G \to \comps$ by letting $\check{f}(x) := f(x^{-1})$ 
for $x \in G$. Let $p \in (1,\infty)$, and let $p' \in (1,\infty)$ be dual to 
$p$, i.e.\ $\frac{1}{p} + \frac{1}{p'} = 1$. The \emph{Fig\`a-Talamanca--Herz 
algebra} $A_p(G)$ consists of those functions $f \!: G \to \comps$
such that there are sequences $( \xi_n )_{n=1}^\infty$ in $L^p(G)$ and 
$( \eta_n )_{n=1}^\infty$ in $L^{p'}(G)$ such that
\begin{equation} \label{Apeq1}
  \sum_{n=1}^\infty \| \xi_n \|_{L^p(G)} \| \eta_n \|_{L^{p'}(G)} < \infty
\end{equation}
and
\begin{equation} \label{Apeq2}
  f = \sum_{n=1}^\infty \xi_n \ast \check{\eta}_n.
\end{equation}
The norm on $A_p(G)$ is defined as the infimum over all sums (\ref{Apeq1}) 
such that (\ref{Apeq2}) holds. It is clear that $A_p(G)$ is a Banach space 
that embeds contractively into ${\cal C}_0(G)$, the algebra of all continuous 
functions on $G$ vanishing at infinity. It was shown by C.\ Herz (\cite{Her1}; see also \cite{EymB} or \cite{Pie}) that $A_p(G)$ is closed under pointwise multiplication and thus a Banach algebra. More specifically (\cite[Proposition 3 and Theorem 3]{Her2}), $A_p(G)$ is a regular, Tauberian, commutative Banach algebra whose character space can be canonically identified with $G$. If $p =2$, the algebra
$A_2(G)$ is Eymard's \emph{Fourier algebra} $A(G)$ (\cite{Eym}). (With our 
notation, we follow \cite{EymB} --- as does \cite{Pie} ---
rather than \cite{Her1} and \cite{Her2} like most authors do: 
$A_p(G)$ in our sense is $A_{p'}(G)$ in \cite{Her1} and \cite{Her2}).
\par
The algebras $A_p(G)$ are related to certain isometric representations of
$G$. Let $\lambda_{p'} \!: G \to {\cal B}(L^{p'}(G))$ be the \emph{regular 
left representation} of $G$ on $L^{p'}(G)$, i.e.\
\[
  (\lambda_{p'}(x)\xi)(y) = \xi(x^{-1} y) \qquad 
  (x,y \in G, \, \xi \in L^{p'}(G))
\] 
The algebra of \emph{$p'$-pseudomeasures} $\PM_{p'}(G)$ is defined as the $w^\ast$-closure of $\lambda_{p'}(L^1(G))$ in the dual Banach space 
${\cal B}(L^{p'}(G))$. There is a canonical duality $\PM_{p'}(G) \cong A_p(G)^\ast$ via
\[
  \langle \xi \ast \check{\eta}, T \rangle := \langle T\eta, \xi \rangle \qquad (\xi \in L^{p'}(G), \, \eta \in L^p(G), \, T \in \PM_{p'}(G)). 
\]
In particular, we have
\[
  (\xi \ast \check{\eta})(x) = \langle \lambda_{p'}(x) \eta, \xi \rangle
  \qquad (\xi \in L^{p'}(G), \, \eta \in L^p(G), \, x \in G). 
\]
Hence, even though it seems to be still unknown (see \cite[9.2]{EymB}) if
$A_p(G)$ consists of coefficient functions of $\lambda_{p'}$ --- except
if $p =2$, of course ---, the elements of $A_p(G)$ are nevertheless not far
from  being coefficient functions of $\lambda_{p'}$ and are, in fact,
coefficient functions of a representation closely related to $\lambda_{p'}$:
\begin{proposition} \label{coeffprop}
Let $G$ be a locally compact group, let $p \in (1, \infty)$, and let
$\lambda_{p'}^\infty \!: G \to {\cal B}(\ell^{p'}(L^{p'}(G)))$ be defined
by letting 
\[
  \lambda_{p'}^\infty(x) := \id_{\ell^p} \tensor \lambda_{p'}(x)
  \qquad (x \in G).
\]
Then $(\lambda_{p'}^\infty, \ell^{p'}(L^{p'}(G)))$ is an isometric 
representation of $G$ and every $f \in A_p(G)$ is a coefficient function of
$(\lambda_{p'}^\infty, \ell^{p'}(L^{p'}(G))$. More precisely, for every 
$\epsilon > 0$, there are $\eta \in \ell^{p'}(L^{p'}(G)))$ and 
$\xi \in \ell^p(L^p(G))$ such that $\| \eta \| \| \xi \| < \| f \| +
\epsilon$ and
\[
  f(x) = \langle \lambda_{p'}^\infty(x) \eta, \xi \rangle \qquad (x \in G).
\]
\end{proposition}
\begin{proof}
To check that $(\lambda_{p'}^\infty, \ell^{p'}(L^{p'}(G)))$ is an isometric 
representation of $G$ is straightforward.
\par
Let $f \in A_p(G)$ and let $\epsilon > 0$. By the definition of $A_p(G)$,
there are sequences $\left( \tilde{\xi}_n \right)_{n=1}^\infty$ in $L^p(G)$ and
$\left( \tilde{\eta}_n \right)_{n=1}^\infty$ in $L^{p'}(G)$ such that
\[
  f(x) = \sum_{n=1}^\infty \left\langle \lambda_{p'}(x)\tilde{\eta}_n,
                           \tilde{\xi}_n \right\rangle \qquad (x \in G)
\]
and
\[
  \sum_{n=1}^\infty \left\| \tilde{\xi}_n \right\|  
                    \left\| \tilde{\eta}_n \right\| < \| f \| + \epsilon.
\]
For $n \in \posints$, set
\[
  \xi_n := \left\{ \begin{array}{cl}
  \left\| \tilde{\xi}_n \right\|^{-1+ \frac{1}{p}} 
  \left\| \tilde{\eta}_n \right\|^\frac{1}{p} \tilde{\xi}_n, 
  & \text{if $\tilde{\xi}_n
  \neq 0$}, \\
  0, & \text{otherwise}
  \end{array} \right.
\]
and
\[
   \eta_n := \left\{ \begin{array}{cl}
  \left\| \tilde{\eta}_n \right\|^{-1+ \frac{1}{p'}}
  \left\| \tilde{\xi}_n \right\|^\frac{1}{p'} \tilde{\eta}_n, 
  & \text{if $\tilde{\eta}_n
  \neq 0$}, \\
  0, & \text{otherwise}.
  \end{array} \right.
\]
It follows that,
\[
  \left( \sum_{n=1}^\infty \| \xi_n \|^p \right)^\frac{1}{p}
  = \left( \sum_{n=1}^\infty \left\| \tilde{\xi}_n \right\| 
  \left\| \tilde{\eta}_n \right\| \right)^\frac{1}{p} 
  < ( \| f \| + \epsilon)^\frac{1}{p}
\]
and, similarly,
\[
  \left( \sum_{n=1}^\infty \| \eta_n \|^{p'} \right)^\frac{1}{p'}
  < ( \| f \| + \epsilon)^\frac{1}{p'}.
\]
Consequently, $\xi := ( \xi_n )_{n=1}^\infty \in \ell^p(L^p(G))$ and
$\eta := ( \eta_n )_{n=1}^\infty \in \ell^{p'}(L^{p'}(G))$ satisfy
\[
  \| \xi \| \| \eta \| < \| f \| + \epsilon
\]
as well as
\[
  f(x) = \sum_{n=1}^\infty \left\langle \lambda_{p'}(x)\tilde{\eta}_n,
                           \tilde{\xi}_n \right\rangle 
  = \sum_{n=1}^\infty \langle \lambda_{p'}(x)\eta_n,
                      \xi_n \rangle 
  = \langle \lambda_{p'}^\infty(x) \eta, \xi \rangle \qquad
  (x \in G).
\]
This completes the proof.
\end{proof}
\begin{remark}
The proof of Proposition \ref{coeffprop} is patterned after that of
\cite[Proposition 5]{Daw}.
\end{remark}
\subsection{Ideals with a bounded approximate identity}
In this section, we shall characterize --- for arbitrary $G$ and  
$p \in (1,\infty)$ --- those closed ideals of $A_p(G)$ that have an
approximate identity bounded by $1$. 
\par
Given a locally compact group $G$, $p \in (1,\infty)$, and a closed subset
$F$ of $G$, we let
\[
  I(F) := \{ f \in A_p(G) : f|_F \equiv 0 \}.
\]
\begin{theorem} \label{BAI}
Let $G$ be a locally compact group, and let $p \in (1,\infty)$.
Then the following are equivalent for a closed ideal $I$ of $A_p(G)$:
\begin{items}
\item $I$ has an approximate identity bounded by $1$;
\item there are $x \in G$ and an open, amenable subgroup $H$ of $G$ such that
$I = I(G \setminus xH)$.
\end{items}
\end{theorem}
\par
Our key to proving Theorem \ref{BAI} is the following proposition:
\begin{proposition} \label{coeffs1}
Let $G$ be a locally compact group, let $p \in (1,\infty)$, and let
$( f_\alpha )_{\alpha \in \mathbb A}$ be a bounded net in $A_p(G)$ that
converges pointwise on $G$ to a function $f \!: G \to \comps$. Then there
is a measure space $X$ and an isometric representation $(\pi, L^{p'}(X))$
of $G_d$ such that $f$ is a coefficient function of $(\pi, L^{p'}(X))$. More
precisely, if $C \geq 0$ is such that $\sup_\alpha \| f_\alpha \| \leq C$,
then there are $\eta \in L^{p'}(X)$ and $\xi \in L^p(X)$ with
$\| \eta \| \| \xi \| \leq C$ and
\[
  f(x) = \langle \pi(x)\eta, \xi \rangle \qquad (x \in G).
\]
\end{proposition}
\par
Before we prove Proposition \ref{coeffs1}, we recall a few facts about
ultrapowers of Banach spaces (see \cite{Hei} and \cite{Sim}).
\par
Let $E$ be a Banach space, and let $\mathbb I$ be any index set. We denote
the Banach space of all bounded families $( \xi_i )_{i \in \mathbb I}$ in $E$,
equipped with the supremum norm, by $\ell^\infty(\mathbb{I},E)$. Let
$\mathfrak U$ be an ultrafilter on $\mathbb I$, and define
\[
  {\cal N}_{\mathfrak U} := \left\{ (\xi_i)_{i \in \mathbb I} 
  \in \ell^\infty(\mathbb{I},E): 
  \lim_{\mathfrak U} \| \xi_i \| = 0 \right\}.
\]
Then ${\cal N}_{\mathfrak U}$ is a closed subspace of 
$\ell^\infty(\mathbb{I},E)$. The quotient space $\ell^\infty(\mathbb{I},E) /
{\cal N}_{\mathfrak U}$ is called an \emph{ultrapower} of $E$ and denoted by
$(E)_{\mathfrak U}$. For any $(\xi_i)_{i \in \mathbb I} \in
\ell^\infty(\mathbb{I},E)$, we denote its equivalence class in 
$(E)_{\mathfrak U}$ by $(\xi_i)_{\mathfrak U}$; it is easy to see that
\begin{equation} \label{ulnorm}
  \| (\xi_i)_{\mathfrak U} \|_{(E)_{\mathfrak U}} = \lim_{\mathfrak U}
  \| \xi_i \|_E.
\end{equation}
\par
We require the following facts about ultrapowers:
\begin{itemize}
\item If $E = L^p(X)$ for $p \in (1,\infty)$ and some measure space $X$, then
  there is a measure space $Y$ such that $(E)_{\mathfrak U} \cong L^p(Y)$
  (\cite[Theorem 3.3(ii)]{Hei}).
\item There is a canonical isometric embedding of $(E^\ast)_{\mathfrak U}$
  into $(E)_{\mathfrak U}^\ast$, via the duality
\[
  \left\langle (\xi_i)_{\mathfrak U} , (\phi_i)_{\mathfrak U} \right\rangle :=
  \lim_{\mathfrak U} \langle \xi_i,\phi_i \rangle
  \qquad
  \left( (\xi_i)_{\mathfrak U} \in (E)_{\mathfrak U}, \,
  (\phi_i)_{\mathfrak U} \in (E^\ast)_{\mathfrak U}
  \right),
\]
which, in general, need not be surjective (\cite[p.\ 87]{Hei}).
\item If $E$ is uniformly convex, then so is $(E)_{\mathfrak U}$ 
(\cite[\S 10, Proposition 6]{Sim}).
\end{itemize}
\begin{proof}[Proof of Proposition \emph{\ref{coeffs1}}]
Let $C \geq 0$ such that $\sup_\alpha \| f_\alpha \| \leq C$.
\par
For each $\alpha \in \mathbb A$ and $\epsilon > 0$, Proposition
\ref{coeffprop} provides $\eta_{\alpha,\epsilon} \in \ell^{p'}(L^{p'}(G))$
and $\xi_{\alpha,\epsilon} \in \ell^p(L^p(G))$ such that 
\begin{equation} \label{coeffestim}
  \| \eta_{\alpha,\epsilon} \| \| \xi_{\alpha,\epsilon} \| \leq C + \epsilon
\end{equation}
and
\begin{equation} \label{coeffeq1}
  f_\alpha(x) = \langle \lambda_{p'}^\infty(x) \eta_{\alpha,\epsilon},
  \xi_{\alpha,\epsilon} \rangle \qquad (x \in G).
\end{equation}
\par
Turn $\mathbb{I} := \mathbb{A} \times (0,\infty)$ into a directed set via
\[
  (\alpha_1,\epsilon_1) \preceq (\alpha_2, \epsilon_2) \defiff
  \text{$\alpha_1 \preceq \alpha_2$ and $\epsilon_1 \geq \epsilon_2$},
\]
and let $\mathfrak U$ be an ultrafilter on $\mathbb I$ that dominates the
order filter. Since $\ell^{p'}(L^{p'}(G)) \cong L^{p'}(\posints \times G)$
is an $L^{p'}$-space, there is a measure space $X$ such that
$(\ell^{p'}(L^{p'}(G)))_{\mathfrak U} \cong L^{p'}(X)$. Define
$\pi \!: G \to {\cal B}(L^{p'}(X))$ by letting
\[
  \pi(x)(\zeta_i)_{\mathfrak U} := ( \lambda_{p'}^\infty(x)
  \zeta_i)_{\mathfrak U} 
  \qquad (x \in G, \, ( \zeta_i )_{\mathfrak U} \in L^{p'}(X)).
\]
It is then clear that $(\pi,L^{p'}(X))$ is an isometric representation of, if
not of $G$, but at least of $G_d$. Set $\eta := 
( \eta_{\alpha,\epsilon} )_{\mathfrak U}$ and $\xi := ( \xi_{\alpha,\epsilon}
)_{\mathfrak U}$, so that $\eta \in L^{p'}(X)$ and $\xi \in L^{p'}(X)^\ast
\cong L^p(X)$. From (\ref{coeffestim}), it is immediate that $\| \eta \| \|
\xi \| \leq C$, and from (\ref{coeffeq1}), we obtain
\[
  f(x) = \lim_{\mathfrak U} f_\alpha(x) = 
  \lim_{\mathfrak U} \langle \lambda_{p'}^\infty(x) \eta_{\alpha,\epsilon},
  \xi_{\alpha,\epsilon} \rangle  
  = \langle \pi(x) \eta, \xi \rangle \qquad (x \in G).
\]
This completes the proof.
\end{proof}
\begin{remark}
The idea to use ultrapowers to ``glue together'' representations of groups
or algebras seems to appear for the first time in \cite{CF} and also ---
less explicitly and, as it seems, independently of \cite{CF} --- in \cite{Daw}.
\end{remark}
\par
Another ingredient of the proof of Theorem \ref{BAI} is:
\begin{lemma} \label{subgroup}
Let $G$ a locally compact group, let $p \in (1,\infty)$, and let 
$H$ be an open subgroup of $G$. Then we have a canonical isometric
isomorphism of $A_p(H)$ and $\{ f \in A_p(G) : \supp(f) \subset H \}$.
\end{lemma}
\begin{proof}
By \cite[Theorem 1]{Her2}, restriction to $H$ is a quotient map
from $A_p(G)$ onto $A_p(H)$. Consequently, we have a contractive inclusion
\[
  \{ f \in A_p(G) : \supp(f) \subset H \} \subset A_p(H).
\]
(This does not require $H$ to be open.)
\par
Since $H$ is open, we may view $L^p(H)$ and $L^{p'}(H)$, respectively,
as closed subspaces of $L^p(G)$ and $L^{p'}(G)$, respectively. From the
definition of $A_p(G)$ and $A_p(H)$ it is then immediate that $A_p(H)$
contractively embeds into $A_p(G)$. 
\end{proof}
\begin{proof}[Proof of Theorem \emph{\ref{BAI}}] (i) $\Longrightarrow$ (ii): Let $F \subset G$ be the \emph{hull} of $I$, i.e.\
\[
  F := \{ x \in G : \text{$f(x) = 0$ for all $f \in I$} \}.
\]
Then $F$ is obviously closed, and $I \subset I(F)$ holds. 
\par
Let $( e_\alpha )_\alpha$ be an approximate identity for $I$ bounded by $1$. 
Let $x \in G \setminus F$. Then there is $f \in I$ such that $f(x) \neq 0$.
Since $\lim_\alpha e_\alpha(x) f(x) = f(x)$, it follows that $\lim_\alpha
e_\alpha(x) = 1$. We conclude that $e_\alpha \to \chi_{G \setminus F}$ 
pointwise on $G$. By Proposition \ref{coeffs1}, there is thus a measure space
$X$ and an isometric representation $(\pi, L^{p'}(G))$ of $G_d$ such that
$\chi_{G \setminus F}$ is a normalized coefficient function of 
$(\pi, L^{p'}(G))$. Since $L^{p'}(X)$ is smooth, Theorem \ref{monico} 
yields that $G \setminus F = xH$ with for some $x \in G$ and a subgroup 
$H$ of $G$. Since $F$ is closed, $xH$ --- and thus $H$ --- must be open.
\par
What remains to be shown is the amenability of $H$. Without loss of
generality, suppose that $F = G \setminus H$, so that
\[
  I \subset I(F) = \{ f \in A_p(G) : \supp(f) \subset H \} \cong A_p(H)
\]
by Lemma \ref{subgroup}. Since the Banach algebra $A_p(H)$ is Tauberian,
and since the hull of $I$ in $H$ is empty, it follows that $I = A_p(H)$,
so that $A_p(H)$ has a bounded approximate identity. By 
\cite[Theorem 6]{Her2}, this means that $H$ is amenable.
\par
(ii) $\Longrightarrow$ (i): Without loss of generality, suppose that $I
= I(G \setminus H)$ for some open subgroup of $G$, so that --- again by
Lemma \ref{subgroup} ---
\[
  I = \{ f \in A_p(G) : \supp(f) \subset H \} \cong A_p(H).
\]
Since $H$ is amenable, $A_p(H)$ has an approximate identity bounded by $1$ 
(\cite[Theorem 6]{Her2}), which proves the claim.
\end{proof}
\begin{remarks}
\item In the $p=2$ case, Theorem \ref{BAI} is \cite[Proposition 3.12]{For1}.
\item For a locally compact group $G$, let
\[
  {\cal R}_c(G) := \{ F \in {\cal R}(G_d) : \text{$F$ is closed} \}.
\]
If $G$ is amenable, then a closed ideal $I$ of
$A(G)$ has a bounded approximate identity if and only if $I = I(F)$ for some
$F \in {\cal R}_c(G)$ (\cite[Theorem 2.3]{FKLS}). The ``if'' part of
this result remains true with $A(G)$ replaced by $A_p(G)$ for arbitrary
$p \in (1,\infty)$ (\cite[Theorem 4.3]{FKLS}). If the Cohen--Host idempotent
theorem could be extended to isometric representations on $L^p$-spaces
for general $p \in (1,\infty)$, then the proof of Theorem \ref{BAI} can
easily be adapted to extend both directions of \cite[Theorem 2.3]{FKLS} to
arbitrary Fig\`a-Talamanca--Herz algebras. 
\end{remarks}
\subsection{Amenability}
The theory of amenable Banach algebras begins with B.\ E.\ Johnson's memoir
\cite{Joh1}. The choice of terminology is motivated by 
\cite[Theorem 2.5]{Joh1}: a locally compact group is amenable (in the usual
sense; see \cite{Pie}, for example), if and only if
its group algebra $L^1(G)$ is an amenable Banach algebra.
\par
Johnson's original definition of an amenable Banach algebra was in terms of
cohomology groups (\cite{Joh1}). We prefer to give another approach, which is 
based on a characterization of amenable Banach algebras from \cite{Joh2}.
\par
Let $\A$ be a Banach algebra, and let $\Tensor$ stand for the (completed) 
projective tensor product of Banach spaces. The space $\A \Tensor \A$ is
a Banach $\A$-bimodule in a canonical manner via\[
  a \cdot (x \tensor y) := ax \tensor y \quad\text{and}\quad (x \tensor y) \cdot a := x \tensor ya \qquad (a,x,y \in \A),
\]
and the \emph{diagonal operator} 
\[
  \Delta_\A \!: \A \Tensor \A \to \A, \quad a \tensor b \mapsto ab
\]
is a homomorphism of Banach $\A$-bimodules.
\begin{definition} \label{diagdef}
Let $\A$ be a Banach algebra. An \emph{approximate diagonal} for $\A$ is
a bounded net $( d_\alpha )_\alpha$ in $\A \Tensor \A$ such that
\begin{equation} \label{diag1}
  a \cdot d_\alpha - d_\alpha \cdot a \to 0 \qquad (a \in \A)
\end{equation}
and
\begin{equation} \label{diag2}
  a \Delta_\A d_\alpha \to a \qquad (a \in \A).
\end{equation}
\end{definition}
\begin{definition} \label{bamdef}
A Banach algebra $\A$ is said to be \emph{$C$-amenable} with $C \geq 1$
if there is an approximate diagonal for $\A$ bounded by $C$. If $\A$ is
$C$-amenable for some $C \geq 1$, then $\A$ is called \emph{amenable}.
\end{definition}
\begin{examples}
\item Let $G$ be a locally compact group. As mentioned already, $L^1(G)$
is amenable in the sense of Definition \ref{bamdef} if and only if $G$ is
amenable, and by \cite{Sto}, $L^1(G)$ is $1$-amenable if and only if $G$ is 
amenable. Hence, for $L^1(G)$, amenability and $1$-amenability are equivalent.
\item A $\cstar$-algebra $\A$ is amenable if and only if it is nuclear
(see \cite[Chapter 6]{LoA} for a relatively 
self-contained exposition of this deep result). In fact, the nuclearity of 
$\A$ implies already that it is $1$-amenable (\cite[Theorem 3.1]{Haa}).
Hence, amenability and $1$-amenability are also equivalent for 
$\cstar$-algebras.
\item In general, $1$-amenability is far more restrictive than mere 
amenability: $A(G)$ is amenable for any finite group $G$, but is $1$-amenable
only if and only if $G$ is abelian (\cite[Proposition 4.3]{Joh3}).
\end{examples}
\par
For more examples and a modern account of the theory of amenable
Banach algebras, see \cite{LoA}.
\par
It is straightforward from Definitions \ref{bamdef} and \ref{diagdef} that
an amenable Banach algebra must have a bounded approximate identity. It
is therefore immediate from Leptin's theorem (\cite{Lep}) and its
generalization to Fig\`a-Talamanca--Herz algebras by Herz 
(\cite[Theorem 6]{Her2}) that, 
for a locally compact group $G$, the Fourier algebra $A(G)$ --- or, more 
generally, $A_p(G)$ for any $p \in (1,\infty)$ --- can be amenable only if
$G$ is amenable. The tempting conjecture that $A(G)$ is amenable if and only
if $G$ is amenable, turned out to be wrong, however: in \cite{Joh3}, Johnson
exhibited examples of compact groups $G$ such that $A(G)$ is not amenable.
Eventually, Forrest and the author (\cite[Theorem 2.3]{FR}) 
gave a characterization of
those $G$ for which $A(G)$ is amenable: they are precisely the \emph{almost
abelian} group, i.e.\ those with an abelian subgroup of finite index.
\par
In \cite{RunPP}, the author gave a more direct proof of \cite[Theorem 2.3]{FR}
that made direct appeal to the Cohen--Host idempotent theorem. Invoking
\cite[Theorem 2.1]{IS}, the author also proved that $A(G)$ is $1$-amenable
for a locally compact group $G$ if and only if $G$ is abelian 
(\cite[Theorem 3.5]{RunPP}). In this section, we shall extend this latter
result to general Fig\`a-Talamanca--Herz algebras.
\par
Even though our arguments in this section, parallel those in the last one,
we now have to consider representations on spaces 
more general than 
$L^p$-spaces (which, nevertheless, will still turn out to be uniformly
convex and have uniformly convex duals). Given two locally compact groups $G$ and $H$, and $p,q \in
(1,\infty)$, the left regular representation of $G \times H$ on $L^p(G,L^q(H))$ is defined as
\[
  \lambda_{p,q} \!: G \times H \to {\cal B}(L^p(G,L^q(H))), \quad
  (x,y) \mapsto \lambda_p(x) \tensor \lambda_q(x).
\]
It is immediate that $(\lambda_{p,q},L^p(G,L^q(H)))$ is an isometric
representation of $G \times H$, as is $(\lambda_{p,q}^\infty,
\ell^p(L^p(G,L^q(H))))$, where
\[
  \lambda_{p,q}^\infty(x,y) := \id_{\ell^p} \tensor \lambda_{p,q}(x,y)
  \qquad (x \in G, \, y \in H).
\]
\par
In analogy with Proposition \ref{coeffprop}, we have:
\begin{lemma} \label{amlem1}
Let $G$ and $H$ be locally compact groups, let $p ,q \in (1, \infty)$, and let
$f \in A_p(G) \Tensor A_q(H)$. Then the Gelfand transform of $f$ on $G \times
H$ is a coefficient function of
$(\lambda_{p,q}^\infty, \ell^p(L^p(G,L^q(H))))$, and for each $\epsilon > 0$, 
there are $\eta \in \ell^{p'}(L^{p'}(G,^{q'}(H)))$ and 
$\xi \in \ell^p(L^p(G,L^q(H)))$ such that 
\begin{equation} \label{coeffineq2}
  \| \eta \| \| \xi \| < \| f \| + \epsilon
\end{equation} 
and
\begin{equation} \label{coeffeq2}
  f(x,y) = 
 \langle \lambda_{p',q'}^\infty(x,y) \eta, \xi \rangle 
  \qquad (x \in G, \, y \in H).
\end{equation}
\end{lemma}
\begin{proof}
Let $\epsilon > 0$.
From the definition of $A_p(G)$ and $A_q(H)$ and from the fact the projective
tensor product is compatible with quotients, it follows that there are 
sequences $( \xi_{n,p} )_{n=1}^\infty$ in $L^p(G)$, $( \xi_{n,q}
)_{n=1}^\infty$ in $L^q(H)$, $(\eta_{n,p} )_{n=1}^\infty$ in $L^{p'}(G)$,
and $( \eta_{n,q} )_{n=1}^\infty$ in $L^{q'}(H))$ such that
\begin{equation} \label{Apqdef1}
  f = \sum_{n=1}^\infty (\xi_{n,p} \ast \check{\eta}_{n,p}) \tensor
  (\xi_{n,q} \ast \check{\eta}_{n,q})
\end{equation}
and
\begin{equation} \label{Apqdef2}
  f = \sum_{n=1}^\infty \| \xi_{n,p} \|_{L^p(G)} 
  \| \check{\eta}_{n,p} \|_{L^{p'}(G)} \| \xi_{n,q} \|_{L^q(H)}
  \| \check{\eta}_{n,q} \|_{L^{q'}(H)} < \| f \| + \epsilon.
\end{equation}
For $n \in \posints$, set
\[
  \xi_n :=  \xi_{n,p} \tensor \xi_{n,q} \in L^p(G,L^q(H))
  \qquad\text{and}\qquad
  \eta_n :=  \eta_{n,p} \tensor \eta_{n,q} \in L^{p'}(G,L^{q'}(H)).
\]
From (\ref{Apqdef1}) and (\ref{Apqdef2}), it is then obvious that
\[
  f(x,y) = \sum_{n=1}^\infty \langle \lambda_{p',q'}(x,y) \eta_n,\xi_n \rangle
  \qquad (x \in G, \, y \in H)
\]
and
\[
  \sum_{n=1}^\infty \| \xi_n \|_{L^p(G,L^q(H))} 
  \| \eta_n \|_{L^{p'}(G,L^{q'}(H))} < \| f \| + \epsilon.
\]
\par
As in the proof of Proposition \ref{coeffprop}, we eventually obtain $\eta \in
\ell^{p'}(L^{p'}(G,L^{q'}(H)))$ and $\xi \in \ell^p(L^p(H,L^q(H)))$ that
satisfy (\ref{coeffineq2}) and (\ref{coeffeq2}).
\end{proof}
\begin{proposition} \label{coeffs2}
Let $G$ and $H$ be locally compact groups, let $p,q \in (1,\infty)$, and let
$( f_\alpha )_{\alpha \in \mathbb A}$ be a bounded net in $A_p(G) \Tensor 
A_q(H)$ that converges pointwise on $G \times H$ to a function $f \!: G 
\times H \to \comps$. Then there is an isometric representation $(\pi, E)$
of $G_d \times H_d$ on a uniformly convex Banach space such that, 
if $C \geq 0$ is such that 
$\sup_\alpha \| f_\alpha \| \leq C$, there are 
$\eta \in E$ and $\xi \in E^\ast$ with $\| \eta \| \| \xi \| \leq C$ and
\[
  f(x,y) = \langle \pi(x,y)\eta, \xi \rangle \qquad (x \in G, \, y \in H).
\]
\end{proposition}
\begin{proof}
The proof parallels that of Proposition \ref{coeffs1}, so that we can afford
being somewhat sketchy.
\par
Let $C \geq 0$ such that $\sup_\alpha \| f_\alpha \| \leq C$.
For each $\alpha \in \mathbb A$ and $\epsilon > 0$, Lemma \ref{amlem1}
provides $\eta_{\alpha,\epsilon} \in \ell^{p'}(L^{p'}(G,L^{q'}(H)))$
and $\xi_{\alpha,\epsilon} \in \ell^p(L^p(G,L^q(H)))$ such that 
\[
  \| \eta_{\alpha,\epsilon} \| \| \xi_{\alpha,\epsilon} \| \leq C + \epsilon
\]
and
\[
  f_\alpha(x,y) = \langle \lambda_{p',q'}^\infty(x,y) \eta_{\alpha,\epsilon},
  \xi_{\alpha,\epsilon} \rangle \qquad (x \in G, \, y \in H).
\]
\par
As in the proof of Proposition \ref{coeffs1}, turn 
$\mathbb{I} := \mathbb{A} \times (0,\infty)$ into a directed set, and let
$\mathfrak U$ be an ultrafilter on $\mathbb I$ that dominates the
order filter. Since $\ell^{p'}(L^{p'}(G,L^{q'}(H)))$ is uniformly convex
by \cite{Day}, so is $E := (\ell^{p'}(L^{p'}(G,L^{q'}(H))))_{\mathfrak U}$. 
Define $\pi \!: G \times H \to {\cal B}(E)$ by letting
\[
  \pi(x,y)(\zeta_i)_{\mathfrak U} := ( \lambda_{p',q'}^\infty(x,y)
  \zeta_i)_{\mathfrak U} 
  \qquad 
  (x \in G, \, y \in H, \, ( \zeta_i )_{\mathfrak U} \in E), 
\]
and set $\eta := ( \eta_{\alpha,\epsilon} )_{\mathfrak U}$ and 
$\xi := ( \xi_{\alpha,\epsilon} )_{\mathfrak U}$.
\par
As in the proof of Proposition \ref{coeffs1}, it is seen that $(\pi,E)$,
$\eta$, and $\xi$ have the desired properties.
\end{proof}
\par
We obtain finally:
\begin{theorem} \label{oneam}
Let $G$ be a locally compact group. Then the following are equivalent:
\begin{items}
\item $G$ is abelian;
\item $A_p(G)$ is $1$-amenable for each $p \in (1,\infty)$;
\item $A(G)$ is $1$-amenable;
\item there is $p \in (1,\infty)$ such that $A_p(G)$ is $1$-amenable.
\end{items}
\end{theorem}
\begin{proof}
(i) $\Longrightarrow$ (iii): Suppose that $G$ is abelian. Then $A(G) \cong
L^1(\hat{G})$ is $1$-amenable (\cite{Sto}).
\par
(iii) $\Longrightarrow$ (ii): Suppose that $A(G)$ is $1$-amenable, and let
$p \in (1,\infty)$. Since $G$ must be amenable, \cite[Theorem C]{Her1}
yields that $A(G)$ is contained in $A_p(G)$ such that the inclusion is 
contractive. A glance at the proof of \cite[Proposition 2.3.1]{LoA} shows
that $A_p(G)$ then must be $1$-amenable, too.
\par
(ii) $\Longrightarrow$ (iv) is trivial.
\par
(iv) $\Longrightarrow$ (i): Let $p \in (1,\infty)$ be such that $A_p(G)$ is
$1$-amenable, and let $( d_\alpha )_\alpha$ be an approximate
diagonal for $A_p(G)$, bounded by $1$. Since
\[
  \mbox{}^\vee \!: A_p(G) \to A_{p'}(G), \quad f \mapsto \check{f}
\]
is an isometric isomorphism of Banach algebras, the net 
$( (\id_{A_p(G)} \tensor 
\mbox{\ }^\vee) d_\alpha )_\alpha$, which lies in $A_p(G) \Tensor
A_{p'}(G)$, is also bounded by $1$. From (\ref{diag1}) and (\ref{diag2}),
it is immediate that 
$( (\id_{A_p(G)} \tensor \mbox{\ }^\vee) d_\alpha )_\alpha$ converges 
pointwise on $G \times G$ to $\chi_{\Gamma}$, where
\[
  \Gamma := \{ (x,x^{-1}) : x \in G \}.
\]
By Proposition \ref{coeffs2}, there
is therefore an isometric representation of $G_d \times G_d$ 
on a uniformly convex Banach space such that $\chi_{\Gamma}$ is a normalized
coefficient function of $(\pi,E)$. From Theorem \ref{monico}, 
we conclude that $\Gamma$ is a left coset of a subgroup of $G \times G$.
Since $\Gamma$ contains the identity of $G \times G$,
it follows that $\Gamma$ is, in fact, a subgroup of $G \times G$.
This is possible only if $G$ is abelian.
\end{proof}
\begin{remark}
Let $G$ be a locally compact group, and consider the following statements:
\begin{items}
\item $G$ is almost abelian;
\item $A_p(G)$ is amenable for each $p \in (1,\infty)$;
\item $A(G)$ is amenable;
\item there is $p \in (1,\infty)$ such that $A_p(G)$ is amenable.
\end{items}
It is known that (i) $\Longrightarrow$ (ii) 
(\cite[Corollary 8.4]{RunExpo}), and (ii) $\Longrightarrow$ (iii)
$\Longrightarrow$ (iv) are trivial. We believe, but have been unable to prove,
that (i), (ii), (iii), and (iv) are equivalent. An inspection of the
proof of Theorem \ref{oneam} reveals that the main obstacle in the way of
proving (iv) $\Longrightarrow$ (i) is the lack of a suitable Cohen--Host type
idempotent theorem for isometric representations on Banach spaces of the form arising in Proposition \ref{coeffs2}.
\end{remark}
\renewcommand{\baselinestretch}{1.0}
\renewcommand{\baselinestretch}{1.2}
\dated
\vfill
\begin{tabbing}
{\it Author's address\/}: \= Department of Mathematical and Statistical Sciences \\
\> University of Alberta \\
\> Edmonton, Alberta \\
\> Canada T6G 2G1 \\[\medskipamount]
{\it E-mail\/}: \> {\tt vrunde@ualberta.ca} \\[\medskipamount]
{\it URL\/}: \> {\tt http://www.math.ualberta.ca/$^\sim$runde/} 
\end{tabbing} 

\begin{thebibliography}{B--F--G--M}
%
\begin{small}
%
%
\bibitem[B--F--G--M]{BFGM} \textsc{U.\ Bader}, \textsc{A.\ Furman}, \textsc{T.\ Gelander}, and \textsc{N.\ Monod}, Property (T) and rigidity for actions on Banach spaces. ArXiv: \texttt{math.GR/0506361}.
%
\bibitem[Coh 1]{Coh} {\sc P.\ J.\ Cohen}, On a conjecture of Littlewood and 
idempotent measures. \textit{Amer.\ J.\ Math.}\ \textbf{82} (1960), 191--212.
%
\bibitem[Coh 2]{Coh2} {\sc P.\ J.\ Cohen}, On homomorphisms of group algebras.
  \textit{Amer.\ J.\ Math.}\ \textbf{82} (1960), 213--226. 
%
\bibitem[C--F]{CF} {\sc M.\ Cowling} and {\sc G.\ Fendler}, On representations 
in Banach spaces. \textit{Math.\ Ann.}\  \textbf{266} (1984), 307--315. 
%
%
\bibitem[Daw]{Daw} {\sc M.\ Daws}, Aren regularity of the algebra of operators
on a Banach space. \textit{Bull.\ London Math.\ Soc.}\ \textbf{36} (2004),
493--503.
%
\bibitem[Day]{Day} {\sc M.\ M.\ Day},  Some more uniformly convex spaces.
\textit{Bull.\ Amer.\ Math.\ Soc.}\ \textbf{47} (1941), 504--507.
%
\bibitem[E--R]{ER} {\sc E.\ G.\ Effros} and {\sc Z.-J.\ Ruan}, 
\textit{Operator Spaces}. London Mathematical Society Monographs (New Series)
\textbf{23}, Clarendon Press, 2000.
%
\bibitem[E--S]{ES} {\sc M.\ Enock} and {\sc J.-M.\ Schwartz}, 
\textit{Kac Algebras and Duality of Locally Compact Groups} (with a preface by 
A.\ Connes and a postface by A.\ Ocneanu). Springer-Verlag, 1992.
%
\bibitem[Eym 1]{Eym} {\sc P.\ Eymard}, L'alg\`ebre de Fourier d'un groupe localement compact. \textit{Bull.\ Soc.\ Math.\ France} \textbf{92} (1964), 181--236. 
%
\bibitem[Eym 2]{EymB} {\sc P.\ Eymard}, Alg\`ebres $A_p$ et convoluteurs de 
$L^p$. In: \textit{S\'eminaire Bourbaki, vol.\ 1969/70, Expos\'es 364--381},
pp.\ 55--72. Lecture Notes in Mathematics \textbf{180}. 
Springer Verlag, 1971.
%
\bibitem[Fab \textit{et al.}]{Fab} {\sc M.\ Fabian}, {\sc P.\ Habala},
{\sc P.\ H\'ajek}, {\sc V.\ Montesinos Santaluc\'{\i}a}, {\sc J.\ Pelant},
and {\sc V.\ Zizler}, \textit{Functional Analysis and Infinite-Dimensional 
Geometry}. CMS Books in Mathematics \textbf{8}, Springer Verlag,
2001. 
%
\bibitem[Fen]{Fen} {\sc G.\ Fendler}, A uniformly bounded representation 
associated to a free set in a discrete group. \textit{Colloq.\ Math.}\
\textbf{59} (1990), 223--229.
%
\bibitem[For 1]{For1} {\sc B.\ E.\ Forrest}, Amenability and bounded
approximate identities in ideals of $A(G)$. {\it Illinois J.\ Math.\/}\
{\bf 34\/} (1990), 1--25.
%
\bibitem[For 2]{For2} {\sc B.\ E.\ Forrest}, Amenability and the structure of 
the algebras $A_p(G)$. \textit{Trans.\ Amer.\ Math.\ Soc.}\ \textbf{343} 
(1994), 233--243.
%
\bibitem[F--K--L--S]{FKLS} {\sc B.\ E.\ Forrest}, {\sc E.\ Kaniuth}, 
{\sc A.\ T.-M.\ Lau}, and {\sc N.\ Spronk}, Ideals with bounded approximate
identities in Fourier algebras. \textit{J.\ Funct.\ Anal.}\ \textbf{203}
(2003), 286--304.
%
\bibitem[F--R]{FR} {\sc B.\ E.\ Forrest} and {\sc V.\ Runde}, Amenability and 
weak amenability of the Fourier algebra. \textit{Math.\ Z.}\ \textbf{250} (2005), 731--744,
%
\bibitem[Gra]{Gra} {\sc E.\ E.\ Granirer}, On some properties of the Banach 
algebras $A_p(G)$ for locally compact groups. \textit{Proc.\ Amer.\ Math.\ 
Soc.}\ \textbf{95} (1985), 375--381. 
%
\bibitem[Haa]{Haa} {\sc U.\ Haagerup}, All nuclear $\cstar$-algebras are 
amenable. {\it Invent.\ math.\/}\ {\bf 74\/} (1983), 305--319.
%
\bibitem[Hei]{Hei} {\sc S.\ Heinrich}, Ultraproducts in Banach space theory.  
\textit{J.\ reine angew.\ Math.}\ \textbf{313} (1980), 72--104.
%
\bibitem[Her 1]{Her1} {\sc C.\ Herz}, The theory of $p$-spaces with an application to convolution operators. \textit{Trans.\ Amer.\ Math.\ Soc.}\ \textbf{154} 
(1971), 69--82.
%
\bibitem[Her 2]{Her2} {\sc C.\ Herz}, Harmonic synthesis for subgroups. 
\textit{Ann.\ Inst.\ Fourier (Grenoble)} \textbf{23} (1973), 91--123. 
%
\bibitem[Hos]{Hos} {\sc B.\ Host}, Le th\'eor\`eme des idempotents dans
$B(G)$. \textit{Bull.\ Soc.\ Math.\ France} \textbf{114} (1986), 215--223. 
%
\bibitem[Ili]{Ili} {\sc M.\ Ilie}, On Fourier algebra homomorphisms.  
\textit{J.\ Funct.\ Anal.}\ \textbf{213} (2004), 88--110.
%
\bibitem[I--S]{IS} {\sc M.\ Ilie} and {\sc N.\ Spronk}, Completely bounded 
homomorphisms of the Fourier algebras. {\it J.\ Funct.\ Anal.\/}\ \textbf{225} (2005), 480--499.
%
\bibitem[Joh 1]{Joh1} {\sc B.\ E.\ Johnson}, Cohomology in Banach algebras. 
\textit{Mem.\ Amer.\ Math.\ Soc.}\ \text{127} (1972).
%
\bibitem[Joh 2]{Joh2} {\sc B.\ E.\ Johnson}, Approximate diagonals and cohomology of certain annihilator Banach algebras.
\textit{Amer.\ J.\ Math.}\ \textbf{94} (1972), 685--698. 
%
\bibitem[Joh 3]{Joh3} {\sc B.\ E.\ Johnson}, Non-amenability of the Fourier algebra of a compact group. \textit{J.\ London Math.\ Soc.}\ (2)
\textbf{50} (1994), 361--374.  
%
\bibitem[J--L]{JL} {\sc W.\ B.\ Johnson} and {\sc J.\ Lindenstrauss},
Basic concepts in the geometry of Banach spaces, pp.\ 1--84. In: 
{\sc W.\ B.\ Johnson} and {\sc J.\ Lindenstrauss} (ed.s), 
\textit{Handbook of the Geometry of Banach spaces}, I. North-Holland
Publishing Co., 2001.
%
\bibitem[L--N--R]{LNR} {\sc A.\ Lambert}, {\sc M.\ Neufang}, and {\sc V.\ Runde}, Operator space structure and amenability for Fig\`a-Talamanca--Herz algebras. \textit{J.\ Funct.\ Anal.}\ \textbf{211} (2004), 245--269.
%
\bibitem[Lei 1]{Lei1} {\sc M.\ Leinert}, Faltungsoperatoren auf gewissen 
diskreten Gruppen. \textit{Studia Math.}\ \textbf{52} (1974), 149--158.
%
\bibitem[Lei 2]{Lei2} {\sc M.\ Leinert}, Absch\"atzung von Normen gewisser
Matrizen und eine Anwendung. \textit{Math.\ Ann.}\ \textbf{240} (1979),
13--19.
%
\bibitem[Lep]{Lep} {\sc H.\ Leptin}, Sur l'alg\`ebre de Fourier d'un groupe 
localement compact. \textit{C.\ R.\ Acad.\ Sci.\ Paris}, S\'er.\ A \textbf{266}
(1968), 1180--1182. 
%
\bibitem[Mia]{Mia} {\sc T.\ Miao}, Predual of the multiplier algebra of 
$A_p(G)$ and amenability. \textit{Canad.\ J.\ Math.}\ \textbf{56} (2004),  
344--355. 
%
\bibitem[Pie]{Pie} {\sc J.\ P.\ Pier}, \textit{Amenable Locally Compact 
Groups}. Wiley-Interscience, 1984.
%
\bibitem[Rua]{Rua} {\sc Z.-J.\ Ruan}, The operator amenability of
$A(G)$. \textit{Amer.\ J.\ Math.}\ \textbf{117} (1995), 1449--1474.
%
\bibitem[Run 1]{LoA} {\sc V.\ Runde}, \textit{Lectures on Amenability}. 
Lecture Notes in Mathematics \textbf{1774}, Springer Verlag, 2002.
%
\bibitem[Run 2]{RunExpo} {\sc V.\ Runde}, Applications of operator spaces to 
abstract harmonic analysis. \textit{Expo.\ Math.}\ \textbf{22} (2004), 
317--363.
%
\bibitem[Run 3]{RunPP} {\sc V.\ Runde}, The amenability constant of the
Fourier algebra. {\it Proc.\ Amer.\ Math.\ Soc.\/}\ \textbf{134} (2006), 1473--1481.
%
\bibitem[Sim]{Sim} {\sc B.\ Sims}, \textit{``Ultra''-techniques in Banach space
theory}. Queen's Papers in Pure and Applied Mathematics \textbf{60}, Queen's 
University, 1982.
%
\bibitem[Sto]{Sto} {\sc R.\ Stokke}, Approximate diagonals and F{\o}lner 
conditions for amenable group and semigroup algebras. {\it Studia Math.\/}\
{\bf 164\/} (2004), 139--159.
%
\end{small}
%
\end{thebibliography}
\end{document}